%% file: legendre2.tex
\documentclass[11pt]{amsart}
\usepackage[margin=1in]{geometry}
\usepackage{amsthm}
\usepackage{times,eucal,fullpage,times,amsmath,amsthm,amssymb,mathrsfs,stmaryrd,enumerate,accents}
\usepackage[all]{xy}
\usepackage[usenames,dvipsnames]{color}
\usepackage[colorlinks,citecolor=blue,hyperindex]{hyperref}
\usepackage{epigraph}
\usepackage{graphicx}
\usepackage{fancyhdr} 
\usepackage{enumitem}
\input{shortforms}

\newtheorem{theorem}[equation]{Theorem}      
\newtheorem{lemma}[equation]{Lemma}          %
\newtheorem{corollary}[equation]{Corollary}  
\newtheorem{proposition}[equation]{Proposition}

\theoremstyle{definition}
\newtheorem{conj}[equation]{Conjecture}

\theoremstyle{definition}
\newtheorem{defn}[equation]{Definition}
\theoremstyle{remark}

\theoremstyle{definition}
\newtheorem{remark}[equation]{Remark}

\numberwithin{equation}{section}

\let\into=\hookrightarrow

\newcommand{\Q}{{\mathbb Q}}

\newcommand{\Z}{{\mathbb Z}}

\renewcommand{\int}{\operatorname{int}}
\renewcommand{\O}{{\mathcal O}}

\renewcommand{\wp}{{\mathfrak p}}

\newcommand{\fq}{{\mathbb{F}_q}}
\newcommand{\fqt}{\fq(t)}

\renewcommand{\bpro}{\begin{proposition}}
\renewcommand{\epro}{\end{proposition}}

\long\def\comment#1\endcomment{}

\begin{document}

\title[]{A method for construction of rational points over elliptic curves II\\
	Points over solvable extensions}%
\author{Kirti Joshi}%
\address{Math. department, University of Arizona, 617 N Santa Rita, Tucson
85721-0089, USA.} \email{kirti@math.arizona.edu}

\thanks{}%


\begin{abstract}
I provide a systematic construction of points, defined over finite radical extensions,  on any Legendre curve over any field. This  includes as special case Douglas Ulmer's construction of rational points over a rational  function field in characteristic $p>0$.  In particular I show that if $n\geq 4$ is any even integer and not divisible by the characteristic of the field then any elliptic curve $E$ over this field  has at least $2n$ rational points over a finite solvable field extension.  Under additional hypothesis, I show that these are   of infinite order. I also show that Ulmer's points lift to characteristic zero and in particular to the canonical lifting. 
\end{abstract}
\maketitle
\setlength{\epigraphwidth}{.45\textwidth}
\epigraph{With no-mind, blossoms invite the butterfly;\\
	With no-mind, the butterfly visits the blossoms.\\
	When the flowers bloom, the butterfly comes;\\
	When the butterfly comes, the flowers bloom.}{Ry\={o}kan `Taigu' \cite{ryokan}}
\tableofcontents

\renewcommand{\wp}{\mathfrak{p}}
\newcommand{\wq}{\mathfrak{q}}
\newcommand{\ok}{\O_K}

\section{Introduction}
In \cite{joshi17-legendre} I gave a construction of $n\geq 7$ (odd) rational points on any Legendre elliptic curve over any number field with the $n$ points being defined over large highly non-solvable extensions of the number field. 
In this note I provide a different construction of $n \geq 4$ (even) rational points on any Legendre
elliptic curve over any field of characteristic not equal to two with points being defined over an explicit solvable
extensions. My construction includes as a special case D. Ulmer's classic characteristic $p>0$ construction
(\cite{ulmer14}) and has the additional salient feature of providing explicit rational points on Legendre elliptic
curves over any number field. Under additional hypothesis I show that these points are of infinite order (for
number fields). As far as I am aware no general constructions of this type (and generality) is known for
number fields. Conjecture~\ref{con:linear-rank-conj} is motivated by Ulmer's work.
\section{Points on a Legendre elliptic  curve}
Let $K$ be a field. Throughout this paper I will assume that  $2\neq0$ in $K$. 
A Legendre elliptic over $K$ is the following elliptic curve
\be\label{eq:legendre} 
y^2=x(x+1)(x+\lambda)
\ee
for some $\lambda\in K-\{0,-1\}$. Let me begin with the following basic observation.
\newcommand{\vl}{\sqrt{\frac{u^{n-1}+1}{u+1}}}
\newcommand{\lep}{\left(u,u(u+1)\vl\right)}
\bthm\label{th:main-legendre}
For any even integer $n\geq 4$, and for any Legendre elliptic curve \eqref{eq:legendre} over $K$, the point 
$$\lep$$
is a $K\left[u,\vl\right]$-valued point on \eqref{eq:legendre} with $u^n=\lambda$.
\ethm
This includes as a special case the fundamental construction of D.~Ulmer (see \cite{ulmer14}):
\bthm[D.~Ulmer]\label{th:ulmer}
Let $K$ be a field of characteristic $p\geq 3$ and let $n=p^f+1$ with $f\geq 1$. Then \eqref{eq:legendre} has a $K[u]$-rational point
$$P=\left(u,u(u+1)^{(p^f+1)/2}\right)\text{ where } u^n=u^{p^f+1}=\lambda.$$
\ethm

\bp[Proof of Theorem~\ref{th:main-legendre}]
The assertion is clear from:
\beas
\left(u(u+1)\vl\right)^2&=&u^2(u+1)^2\left(\frac{u^{n-1}+1}{u+1}\right)\\
&=&u(u+1)(u^n+u)\\
&=&u(u+1)(u+\lambda) \text{ as } u^n=\lambda.
\eeas
\ep

\bp[Proof of Theorem~\ref{th:ulmer}]
Noting that $n=p^f+1$ is even as $p$ is odd and $K$ has characteristic $p>0$, the assertion is clear from:
\beas
\left((u+1)\vl\right)&=&(u+1)\sqrt{\frac{u^{p^f}+1}{u+1}}\\
&=&(u+1)\sqrt{\frac{(u+1)^{p^f}}{u+1}}\\
&=&(u+1)\sqrt{(u+1)^{p^f-1}}\\
&=&(u+1)(u+1)^{(p^f-1)/2}\\
&=&(u+1)^{(p^f+1)/2}.	
\eeas
Hence $P=\left(u,u(u+1)\vl \right)=(u,u(u+1)^{(p^f+1)/2})$ is a point defined over $K[u]$.
\ep

\bcor\label{cor:legendre-solvable}
Let $E$ be a Legendre elliptic curve given by \eqref{eq:legendre} over a field $K$ in which $2\neq 0$. Let $n\geq 4$ be any even integer not divisible by the characteristic of $K$. Then there exists a radical extension $K'$ of $K$ over which $E$ has $n$ $K'$-rational points. 
\ecor
\bp 
Let $u_1,\ldots,u_n$ be the $n$ roots of $X^n-\lambda$. Then by Theorem~\ref{th:main-legendre}  the points $\left(u_i,u_i(u_i+1)\sqrt{\frac{u_i^{n-1}+1}{u_i+1}}\right)$ are on $E$ and defined over 
$$
K'=K\left(u_1,\ldots,u_n,\sqrt{\frac{u_1^{n-1}+1}{u_1+1}},\ldots,\sqrt{\frac{u_n^{n-1}+1}{u_n+1}}\right),$$
which is a radical extension of $K$ as $u_i=\zeta^i \sqrt[n]{\lambda}$ where $0\leq i<n$ and $\zeta$ is a primitive $n^{th}$-root of unity.
\ep

The following theorem shows that these considerations can be applied to any elliptic curve (not necessarily a Legendre elliptic curve) over any field $K$ in which $6\neq 0$.

\bthm\label{th:main-solvable-points}
Let $K$ be a field in which $6\neq 0$  and let $E$ be any elliptic curve over $K$. Let $n\geq 4$ be any even integer not divisible by the characteristic of $K$. Then there exists an explicit finite, solvable extension $K'$ of $K$ over which $E$ has least $2n$ rational points.
\ethm
\bp 
Let $E$ be any elliptic curve over $K$. After passage to the field $K(E[2])$ (which is  solvable with Galois group a subgroup of $S_3$) generated by $2$-torsion , one can assume that $E$ has two torsion defined over $K(E[2])$ and therefore $E$ can be defined over this field by $y^2=(x+a)(x+b)(x+c)$ and after applying a suitable automorphism defined over $K(E[2])$ one can assume that $E$ is a Legendre curve \eqref{eq:legendre} over $K(E[2])$. Now the assertion is clear from the construction of the points $P$ carried out over $K(E[2])$ using Theorem~\ref{th:main-legendre} and Corollary~\ref{cor:legendre-solvable}.
\ep

\section{Ulmer's rational points lift to characteristic zero}
My next result is the following consequence of Theorem~\ref{th:main-legendre} which shows that the rational points constructed by D.~Ulmer in \cite{ulmer14} in fact lift to characteristic zero.
\newcommand{\tR}{\tilde{R}}
\newcommand{\sE}{\mathscr{E}}
\newcommand{\tE}{\tilde{E}}
\newcommand{\ttl}{\tilde{t}}
\newcommand{\tP}{\tilde{P}}
\bthm\label{th:ulmer-lifting-1} Let $p$ be an odd prime. Suppose $A/\Z_p$ is a $\Z_p$-algebra of characteristic zero such that $A_0=A/pA\supset \fqt$. Let $\ttl\in A$ be such that $\ttl\cong t\bmod{pA}$. Let $\tE:y^2=x(x+1)(x+\ttl)$ be a Legendre elliptic curve over $\Z_p[\ttl\,]\subset A$. Let $f\geq 1$ be an integer and let $n=p^f+1$. 
Let $B$ be the $A$-algebra given by 
$$B=A\left[u,\sqrt{\frac{u^{n-1}+1}{u+1}}\right] \text{ where } u^n=\ttl.$$
Then
\benum
\item $\tP=\left(u,u(u+1)\sqrt{\frac{u^{n-1}+1}{u+1}}\right)\in \tE(B)$,
\item $\tP\cong (u,u(u+1)^{(p^f+1)/2})\bmod{pB}$ i.e. $\tP\bmod{pB}$ is Ulmer's rational point.
\eenum
\ethm
\bp
The proof is clear from the proofs of Theorem~\ref{th:main-legendre} and Theorem~\ref{th:ulmer}.
\ep

The following Lemma should be well-known but I do not know a reference.
\newcommand{\tl}{\tilde{\lambda}}
\blem 
Let $E$ be a Legendre elliptic curve given by \eqref{eq:legendre} over an algebraically closed field $k$ of characteristic not dividing $6$. Suppose $E$ is ordinary. Let $E^{can}/W(k)$ be the canonical lifting of $E$ to $W(k)$. Then there exists $\tl=\lambda_{can}\in W(k)$,  and $\tl\cong\lambda\bmod{p}$, such that $E^{can}$ is given by
$$y^2=x(x+1)(x+\tl).$$
Moreover $\tl\in W(k)$ is well-defined up to multiplication by  $u^2$,  a unit of $W(k)$, such that $u\cong\pm1\bmod{p}$.
\elem
\bp 
By the theory of canonical liftings, $E$ has a canonical lift $E^{can}$ to $W$ and by well-known arguments, the canonical lift has a Tate-Weierstrass model. As characteristic $p$ of $k$ does not divide $6$ one sees that $E^{can}$ has a model $y^2=f(X)$ with $f(X)\in W(k)$ is monic of degree three. As $f(X)\cong x(x+1)(x+\lambda)\bmod p$ and the factors on the right are coprime so one sees, by Hensel's Lemma, that $f(X)$ has linear factors in $W(k)[x]$ and so $f(X)=(x+a)(x+b)(x+c)$ with $a\cong0\bmod{p}$, etc. Applying an automorphism of $W(k)[x,y]$ one may further assume that $a=0$. So $E^{can}$ has a model $y^2=x(x+b)(x+c)$ with $b\cong-1\bmod{p}$ etc. and as $k$ is algebraically closed so there exists a unit $u\in W(k)$ such that $u^2= b$ and $u^2\cong -1\bmod{p}$ and hence replacing $x,y$ by $xb=xu^2,yu^3$ one sees that $E^{can}$ has a model of the form $y^2=x(x+1)(x+\tl)$ for some $\tl\in W$ and $\tl\cong \lambda\bmod{p}$. The  rest of the assertion is clear.
\ep
I will call such a $\tl=\lambda_{can}$ provided by this lemma a \emph{canonical Serre-Tate-Legendre coordinate} (or simply a \emph{canonical Legendre coordinate}) for  an ordinary, Legendre elliptic curve $E$ given by \eqref{eq:legendre}. The following corollary is now immediate from Theorem~\ref{th:ulmer-lifting-1}:
\newcommand{\tu}{\tilde{u}}
\bthm\label{th:ulmer-lifting-2}
Let $E$ be  the (ordinary) Legendre curve over an algebraic closure $k$ of $\fqt$ with $\lambda=t$ and $E^{can}$ be its canonical lifting with Legendre canonical coordinate $\tl\in W(k)$. Let $n=p^f+1$ for integers $f\geq 1$.   Then Ulmer's rational points $(u,u(u+1)^{p^f+1})$ lifts to a point  $$\tP = \left(\tu,\tu(\tu+1)\sqrt{\frac{\tu^{n-1}+1}{\tu +1}}\right)\in E^{can}(B)$$ where $$B=W(k)\left[\tu,\sqrt{\frac{\tu^{n-1}+1}{\tu +1}}\right] \text{ where } \tu^n=\tl.$$
\ethm

\section{Legendre elliptic curves over number fields}
In this section I assume that $K$ is a number field. Theorem~\ref{th:main-solvable-points} provides a construction of solvable points on any elliptic curve over $K$. Now let me show that under certain circumstances these points are also of infinite order. Analysis is greatly simplified if one assumes that $\lambda$ is an algebraic integer, but it is possible to relax this assumption with additional notational complexity.

\bthm\label{th:main-nf-case}
Let $K$ be a number field and let $E$ given by \eqref{eq:legendre} be a Legendre elliptic curve over $K$ with $\lambda\in\ok$. Suppose $n\geq 4$ is an even integer and let $$L=K\left(u,\vl\right) \text{ with } u^n=\lambda.$$ Suppose $\wp,\wq$ are prime ideals of $\O_L$ not lying over $(2)\subset\Z$ such that
\benum[label={{\bf(\arabic{*})}}]
\item\label{th:main-nf-case-0} $E$ has good reduction at $\wp,\wq$,
\item\label{th:main-nf-case-1} $\vl\in\wp\cap\wq$,
\item\label{th:main-nf-case-2} $u,(u+1)\not\in\wp\cup\wq$.
\eenum
Then $P$ is of infinite order in $E(L)$.
\ethm

\newcommand{\ol}{\O_L}

\bp 
Let $O$ be the point at infinity on $E$. Suppose $P$ is of finite order, say $m$. One can assume that $m\geq 3$ as $P$ is not equal to $O$ nor is it of order two as its $y$-coordinate is non-zero.  Since $E$ has good reduction at $\wp,\wq$ by hypothesis, the equation for $E$ is minimal at $\wp,\wq$. Let $p=\wp\cap\Z$ and $q=\wq\cap \Z$ be the primes of $\Z$ lying below $\wp,\wq$. By assumption $p,q>2$. If $(m,p)=1$ then $E[m]\into E(\ol/\wp)$ (by \cite{silverman-arithmetic}). As $\vl\in\wp$ so $P=(u,0)\bmod{\wp}$ which has order two so $m=2$. Note that as $$u+\lambda=u(u+u^{n-1})=u(u+1)\left(\frac{u^{n-1}+1}{u+1}\right)$$ it follows that $P\cong(-\lambda,0)\bmod{\wp}$, and this also holds  for $\wq$. Thus if $P$ is of finite order $m$ then $p|m$. Now let $m=m'p^r$ where $(m',p)=1$ and $r\geq 1$. Then $Q=m'P$ has order $p^r$. As $p\neq q$ so $E[p^r]\into E(\ol/\wq)$. But again $P\cong (u,0)\bmod{\wq}$ is two torsion,  and reduction modulo $\wq$ is a homomorphism of groups, so modulo $\wq$, $Q=m'P$ is equal to $(u,0)$ if $m'$ is odd or equal to $O$ if $m'$ is even . Thus $Q=Q_0+Q_1$ where $Q_0$ has order  dividing $2$ and $Q_1$ has order equal to a power of $q$ (i.e. $Q_1$ is a torsion element of the kernel of the reduction modulo $\wq$ map). But as $Q$ has order $p^r$ and $Q_0+Q_1$ is  annihilated by $2q^n$ for some $n\geq 0$. But this is clearly  impossible.  Hence one has a contradiction and hence $P$ is not of finite order. So $P$ is of infinite order.
\ep

\bcor 
Suppose $K$ is a number field and $E$ is a Legendre elliptic curve given by \eqref{eq:legendre} for some $\lambda\in\ok$. Assume $n\geq 4$ is an even integer.  Then 
\benum
\item there exists a finite solvable extension $K'\supset K$ such that $E(K')$ has $n$, $K'$-rational points.
\item if conditions \ref{th:main-nf-case-0}--\ref{th:main-nf-case-2} of Theorem~\ref{th:main-nf-case}  hold for every root $u$ of $X^n-\lambda$, then these points are all of infinite order in $E(K')$.
\eenum
\ecor

I do not know if conditions \ref{th:main-nf-case-0}--\ref{th:main-nf-case-2} of Theorem~\ref{th:main-nf-case} hold for a given $\lambda$ and all but finite many even integers $n\geq 4$. But perhaps the following weaker assertion does hold:
\begin{conj}\label{con:linear-rank-conj}
Let $K$ be a number field and assume  that $E$ given by \eqref{eq:legendre} is a Legendre elliptic curve over $K$ with some $\lambda\in K$. Then there exist infinitely many even integers $n\geq 4$ such that conditions \ref{th:main-nf-case-0}--\ref{th:main-nf-case-2} of Theorem~\ref{th:main-nf-case} hold for all the $n$ points $P$ constructed above and these $n$ points generate a subgroup whose rank grows linearly in $n$.
\end{conj}

\section{A numerical example of Theorem~\ref{th:main-nf-case}} The following numerical example shows that there are $\{K,\lambda,n,\wp,\wq\}$ which satisfy all the hypothesis of Theorem~\ref{th:main-nf-case} (so the assertion is non vacuous).

Let $K=\Q$ and consider the curve \eqref{eq:legendre} with $\lambda=86$ and $n=10$. Let $L=\Q(u,v)$ where $$u^n=\lambda=86$$ and $$v^2=\frac{u^{n-1}+1}{u+1}=\frac{u^{9}+1}{u+1}=u^8 - u^7 + u^6 - u^5 + u^4 - u^3 + u^2 - u + 1.$$
Using \cite{sage} one find the factorization into prime ideals:
$$ (v)= (7, v) (37, v) (1069, v) (10934266789, v) (3027381380137219, v).$$
So that choosing $\wp=(37,v)$ and $\wq=(1069,v)$ one sees that $v\in\wp\cap\wq$ and $u,u+1\not\in \wp\cup\wq$. So the hypothesis of Theorem~\ref{th:main-nf-case} are satisfied. In fact one can compute the order of $E(L)_{tor}=8$ ($\lambda$ is a square in $L$ and so $E$ has a four torsion point over $L$) and one checks that $P=(u,u(u+1)v)$ is a non-torsion point in $E(L)$ (as predicted by Theorem~\ref{th:main-nf-case}).

\bibliographystyle{plain}
\bibliography{../points/points.bib,../../master/joshi.bib,../../master/master6.bib}
\end{document}

%% file: shortforms.tex
\newcommand{\be}{\begin{equation}}
\newcommand{\ee}{\end{equation}}
\newcommand{\bes}{\begin{equation*}}
\newcommand{\ees}{\end{equation*}}
\newcommand{\bea}{\begin{eqnarray}}
\newcommand{\eea}{\end{eqnarray}}
\newcommand{\beas}{\begin{eqnarray}}
\newcommand{\eeas}{\end{eqnarray}}
\newcommand{\ben}{\begin{note}}
\newcommand{\een}{\end{note}}
\newcommand{\bexl}{\vskip0.1em\noindent\hrulefill\vskip1em\begin{ExerciseList}}
\newcommand{\eexl}{\end{ExerciseList}\hrulefill}

\newcommand{\bthm}{\begin{theorem}}
\newcommand{\ethm}{\end{theorem}}
\newcommand{\bpro}{\begin{prop}}
\newcommand{\epro}{\end{prop}}
\newcommand{\bcor}{\begin{corollary}}
\newcommand{\ecor}{\end{corollary}}
\newcommand{\bcon}{\begin{conjecture}}
\newcommand{\econ}{\end{conjecture}}
\newcommand{\bp}{\begin{proof}}
\newcommand{\ep}{\end{proof}}
\newcommand{\blem}{\begin{lemma}}
\newcommand{\elem}{\end{lemma}}
\newcommand{\bn}{\begin{note}}
\newcommand{\en}{\end{note}}
\newcommand{\benum}{\begin{enumerate}}
\newcommand{\eenum}{\end{enumerate}}
\newcommand{\bed}{\begin{defn}}
\newcommand{\eed}{\end{defn}}
\newcommand{\brem}{\begin{remark}}
\newcommand{\erem}{\end{remark}}

\newcommand{\btik}{\begin{tikzpicture}\begin{axis}[scale=0.5,axis y line=center, axis x line=middle]}
\newcommand{\etik}{\end{axis}\end{tikzpicture}}

\let\into=\hookrightarrow

\let\cong=\equiv

\newcommand{\upperRomannumeral}[1]{\uppercase\expandafter{\romannumeral#1}}